\documentclass[11pt, a4paper,notitlepage]{article}

\usepackage{amssymb}
\usepackage{amsmath}
\usepackage{amsthm}

\newcommand{\Var}{\mathrm{Var}\,}
\newcommand{\eps}{\varepsilon}

\newcommand{\half}{\frac{1}{2}}

\newtheorem{thm}{Theorem}[section]
\newtheorem{lem}[thm]{Lemma}

\theoremstyle{definition}
\newtheorem{rem}[thm]{Remark}
\newtheorem{ex}[thm]{Example}

\newenvironment{prf}{\begin{proof}[\bf Proof]}{\end{proof}}

\def\infint{\int_{-\infty}^\infty}
\def\ex{{\rm E\,}}
\def\var{{\rm Var\,}}

\def\arg{{\rm Arg\,}}
\newcommand{\cov}{{\rm Cov}\,}

\begin{document}

\title{Nonparametric  Volatility Density Estimation}

\author{Bert van Es, Peter Spreij \\[.5ex]
{\normalsize\sl Korteweg-de Vries Institute for Mathematics}\\
{\normalsize\sl University of Amsterdam}\\
{\normalsize\sl Plantage Muidergracht 24}\\
{\normalsize\sl 1018 TV Amsterdam}\\
{\normalsize\sl The Netherlands}\\[2ex] \and
Harry van Zanten\\[.5ex]
{\normalsize\sl Division of Mathematics and Computer Science}\\
{\normalsize\sl Faculty of Sciences}\\
{\normalsize\sl Free University Amsterdam}\\
{\normalsize\sl De Boelelaan 1081 a}\\
{\normalsize\sl 1081 HV Amsterdam}\\
{\normalsize\sl The Netherlands}}

\maketitle

\begin{abstract}
In this paper we consider a continuous-time stochastic volatility model.
The model contains a stationary volatility process, the
density of which, at a fixed instant in time, we aim to estimate.
We assume that we observe the process at discrete instants in time.
The sampling times will   be
equidistant with vanishing distance.

A Fourier-type deconvolution kernel density estimator based on the
logarithm of the squared processes is proposed to estimate the
volatility density. An expansion of the bias and a bound on the
variance are derived.
\end{abstract}

\section{Introduction}

Let $S$ denote the log price process of some stock in a financial
 market. It is often assumed that $S$ can be modelled as the
solution of a stochastic differential equation or, more general,
as an It\^o diffusion process. So we assume that we can write
\begin{equation}\label{eq:s}
dS_t = b_t\, dt + \sigma_t \,dW_t, \ \ \ S_0=0,
\end{equation}
or, in integral form,
\begin{equation}\label{eq:si}
S_t = \int_0^t b_s\, ds + \int_0^t \sigma_s \,dW_s,
\end{equation}
where $W$ is a standard Brownian motion and the processes $b$ and
$\sigma$ are assumed to satisfy certain regularity conditions
(see Karatzas and Shreve (1991)) to have  the integrals in
(\ref{eq:si}) well-defined. In the financial context, the process
$\sigma$ is called the volatility process.

In this paper we model $\sigma$ as a strictly stationary positive
process satisfying a mixing condition, for example an ergodic
diffusion on $[0,\infty)$, and we make the assumption that
$\sigma$ is independent of $W$. We will assume that the one-dimensional
marginal distribution of $\sigma$  has a density
with respect to the Lebesgue measure on $(0,\infty)$. This is
typically the case in virtually all stochastic volatility models
that are proposed in the literature, where the evolution of $\sigma$ is
modelled by a stochastic differential equation, mostly in terms
of $\sigma^2$, or $\log \sigma^2$ (cf.\ e.g.\ Wiggins (1987), Heston (1993)).

For stochastic
differential equations of the type
\[
d\,X_t = b(X_t)\, dt + a(X_t)\, d\, B_t,
\]
with $B_t$ equal to Brownian motion,
the invariant density is up to a multiplicative constant equal to
\begin{equation}
\label{eq: piet}
x \mapsto \frac{1}{a^2(x)}\,{\exp \left(2\int_{x_0}^x
\frac{b(y)}{a^2(y)}dy\right)},
\end{equation}
where $x_0$ is an arbitrary element of the state space $(l,r)$, see
 e.g.\ Gihman and Skorohod (1972) or Skorokhod (1989). From
formula~(\ref{eq: piet}) one sees that the invariant distribution of the
volatility process (take $X$ for instance equal to $\sigma^2$ or $\log \sigma^2$)
may take  on many different forms,
as is the case for the various models that have been
proposed in the literature.
This observation supports our point of view that nonparametric
procedures are by all means sensible tools to get some insight in
the behaviour of the volatility.

In the present paper we propose a
nonparametric  estimator for the volatility density.
 Using ideas from
deconvolution theory, we will propose a procedure for the
estimation of the marginal density at a fixed point.
We will assume that we observe the log-asset price $S$ at time instants
$0, \Delta, 2\Delta, \ldots, n\Delta$, where the time gap
satisfies $\Delta = \Delta_n \to
0$ and $n\Delta_n \to \infty$ as $n \to \infty$. To asses the
quality of our procedure, we will study how the bias and variance
of the estimator behave under these assumptions.

The remainder of the paper is organized as follows. In the next section,
we give the heuristic arguments that motivate the definition of our estimator.
In Section \ref{sec:results} the main result concerning the asymptotic behaviour
of the estimator is presented and discussed. The proof of the main theorem is
given in the last two sections.

\section{Construction of the estimator}

To motivate the construction of the estimator, we first consider~(\ref{eq:s})
without the drift term, so we
assume to have
\begin{equation*}
  dS_t = \sigma_t \,dW_t, \ \ \ S_0=0.
\end{equation*}
It is  assumed that we observe the process $S$
at the discrete time instants $0$, $\Delta$,
$2\Delta, \ldots, n\Delta$.
For $i=1, 2, \ldots$ we work, as in
Genon-Catalot et al.\ (1998, 1999), with  the normalized
increments
\[
X^{\Delta}_i = \tfrac{1}{\sqrt{\Delta}}(S_{i\Delta} -
S_{(i-1)\Delta}).
\]
For small $\Delta$, we have the  rough approximation
\[
X^{\Delta}_i = \tfrac{1}{\sqrt{\Delta}}
\int_{(i-1)\Delta}^{i\Delta}\sigma_t\,dW_t \approx
\sigma_{(i-1)\Delta}  \tfrac{1}{\sqrt{\Delta}}(W_{i\Delta} -
W_{(i-1)\Delta}) = \sigma_{(i-1)\Delta}Z^\Delta_i,
\]
where for $i=1, 2, \ldots$ we define
\[
Z^\Delta_i= \frac{1}{\sqrt{\Delta}}(W_{i\Delta} -
W_{(i-1)\Delta}).
\]
By the independence and stationarity of Brownian increments, the
sequence $Z^\Delta_1, Z^\Delta_2, \ldots$ is an i.i.d.\ sequence
of standard normal random variables. Moreover,  the sequence is
independent of the process $\sigma$ by assumption.

Taking the logarithm of the square   of
 $X_i^\Delta$  we get
\begin{equation}
\log( (X^\Delta_i)^2)  \approx \log(
\sigma_{(i-1)\Delta}^2) + \log ( (Z_i^{\Delta})^2),
\end{equation}
where the terms in the sum are independent.
Assuming that
the approximation is sufficiently accurate we can use this approximate
convolution structure to
estimate the unknown
density $f$ of $\log(  \sigma_{i\Delta}^2)$ from the observed
$\log ((X^\Delta_i)^2)$.

Before we can define the estimator, we need some more notation.
Observe that the density of the `noise' $\log (Z_i^\Delta)^2$, denoted by $k$, is given by
\begin{equation}\label{densityk}
k(x) = \frac{1}{\sqrt{2\pi}}\, e^{\tfrac{1}{2}x} e^{-
\tfrac{1}{2}e^{x}}.
\end{equation}
The characteristic function of the density $k$ is denoted by $\phi_k$.

We will use a function $w$, satisfying the following condition.
For examples of such kernels see Wand (1998).

\medskip
\noindent{\bf Condition W.} Let $w$ be a real symmetric function
with real valued symmetric characteristic function $\phi_w$ with support
[-1,1].
Assume further
\begin{enumerate}
\item
$\infint |w(u)|du < \infty$
,
$\infint w(u)du=1$
,
$\infint u^2|w(u)|du<\infty$
,
\item
$\phi_w(1-t)=At^\alpha + o(t^\alpha),\quad\mbox{as}\ t\downarrow
0 $ for some $\alpha >0$.
\end{enumerate}

Following a well-known approach in statistical deconvolution theory,
we use a {\em
deconvolution kernel density estimator}, see e.g.\ Section 6.2.4 of Wand and Jones (1995).
Having the characteristic functions $\phi_k$ and $\phi_w$ at our disposal,
choosing a positive {\em bandwidth} $h$, we introduce the kernel function
\begin{equation}\label{fourkernel}
v_h(x)={1\over 2\pi}\infint {{\phi_w(s)}\over\phi_k(s/h)}\
e^{-isx}ds
\end{equation}
and  the density estimator
\begin{equation}\label{eq:fnhc}
f_{nh}(x)={1\over nh}\sum_{j=1}^n v_h\left({{x-\log(
(X^\Delta_j)^2)}\over h}\right).
\end{equation}
One   easily verifies that the function $v_h$, and therefore also
the estimator $f_{nh}$, is real-valued.

\section{Results}\label{sec:results}

To derive the asymptotic behaviour of the estimator, we need a
mixing condition on the process $\sigma$.
For the sake of clarity,
we recall  the basic definitions. For a certain process
$X$ let ${\cal F}_a^b$ be the $\sigma$-algebra of events
generated by the random variables $X_t,\ a \le t \le b$. The
mixing coefficient $\alpha(t)$ is defined by
\begin{equation}\label{eq:alphamix}
\alpha(t)=\sup_{A\in {\cal F}_{-\infty}^0,\ B\in {\cal
F}_t^\infty} |P(A\cap B) -P(A)P(B)|.
\end{equation}
The process $X$ is called {\em  strongly mixing} if
$\alpha(t) \to 0$ as $t\to\infty$.

As we mentioned in the introduction, it is common practice to
model the volatility process $V=\sigma^2$ as the stationary,
ergodic  solution of an SDE of the form
\[
d V_t = b(V_t)\,dt + a(V_t)\,d\,B_t.
\]
The mixing condition that we use in Theorem \ref{contasthm} below is satisfied in this setup.
See for instance Corollary 2.1 of
Genon-Catalot et al.\ (2000).

It is easily verified that for such processes it holds that
$\ex |V_t-V_0| = O(t^{1/2})$,
provided that $b \in L_1(\mu)$ and $a \in L_{2}(\mu)$,
where $\mu$ is the invariant probability measure.
Indeed we have  $\ex |V_t-V_0| \leq \ex \int_0^t |b(V_s)|\,ds + (\ex \int_0^t a^2(V_s)\,
ds)^{1/2} = t||b||_{L_1(\mu)} + \sqrt{t} ||a||_{L_2(\mu)}$.
Although we will not assume explicitly that $\sigma^2$ solves an SDE,
the above observation motivates the following condition.\medskip\\
{\bf Condition $\sigma$.}  We have $\ex|\sigma^2_t -
\sigma_0^2| = O(t^{1/2})$ for $t \to 0$.

\bigskip

The following theorem describes the asymptotic behaviour
of our estimator $f_{nh}$. Note that it  also covers the case where there is a drift
$b_t$ present in equation (\ref{eq:s}).
The condition on the drift
is boundedness of  $\ex b_t^2$. This condition is typically
satisfied in realistic models for the log-returns of a stock,
since $b_t$ is the local rate of return and this will be mostly
bounded itself.

\bigskip

\begin{thm} \label{contasthm}
Assume that $\ex b_t^2$ is
bounded.
Let the process $\sigma$ be strongly mixing with coefficient
$\alpha(t)$ satisfying,  for some $0<q<1$,
$$
\int_0^\infty \alpha(t)^q\,dt<\infty,
$$
and suppose that Condition $\sigma$ holds.
Let the kernel function $w$ satisfy Condition W and let the
density $f$ of\/ $\log\sigma^2_t$ be continuous,
twice continuously differentiable with a bounded second
derivative.
Also assume that the density of $\sigma^2_t$ is bounded in a neighbourhood of zero.
Suppose that
$\Delta=n^{-\delta}$ for  given $0<\delta <1$ and choose
$h=\gamma\pi/\log n$, where $\gamma > 4/\delta$. Then the bias of the estimator
(\ref{eq:fnhc}) satisfies
\begin{equation}\label{contasthm:1}
\ex f_{nh}(x) - f(x) = \tfrac{1}{2}h^2f''(x)\int u^2w(u)du +
o(h^2).
\end{equation}
Moreover, the variance of the estimator satisfies
\begin{equation}\label{contasthm:2}
\var f_{nh}(x) = O\Big({1\over
n}\,h^{2\alpha}e^{\pi/h}\Big)+O\Big({1\over
{nh^{1+q}\Delta}}\Big).
\end{equation}
\end{thm}

\bigskip

The proof of the theorem is deferred to the next section. We conclude the present
section by a number of comments on the result.

\begin{rem} The expectation of the deconvolution estimator is equal to the
expectation of an ordinary kernel density estimator, as becomes clear from the proof of
Lemma~\ref{lem:1}.

It is well-known that the variance of kernel-type deconvolution
estimators heavily depends on the rate of decay to zero of
$|\phi_k(t)|$ as $|t|\to\infty$. The faster the decay the larger
the asymptotic variance. In other words, the smoother $k$ the
harder the estimation problem. This follows for instance for
i.i.d.\ observations from results in Fan (1991) and for stationary
observations from the work of Masry (1993).

The rate of decay of $|\phi_k(t)|$ for the density
(\ref{densityk}) is given by Lemma \ref{phiexpan} in Section
\ref{sec:tech}, stating
$|\phi_k(t)|\sim\sqrt{2}\,e^{-\frac{1}{2}\pi |t|}$, as
$|t|\to\infty$. This shows that $k$ is supersmooth, cf.\ Fan (1991).
By the similarity of the tail of this characteristic function to
the tail of a Cauchy characteristic function  we can expect the
same order of the mean squared error as in Cauchy deconvolution
problems, where it decreases logarithmically in $n$, cf. Fan
(1991) for results on i.i.d.\ observations. Note that this rate,
however slow, is faster than the one for normal deconvolution.
Fan (1991) also shows that we cannot expect anything better.
\end{rem}

\begin{rem}
The choices $\Delta=n^{-\delta}$, with $0<\delta<1$ and
$h=\gamma\pi/\log n$, with $\gamma>4/\delta$ render a variance
that is of order $n^{-1+1/\gamma}(1/\log n)^{2\alpha}$ for the
first term of (\ref{contasthm:2}) and $n^{-1+\delta}(\log
n)^{1+q}$ for the second term. Since by assumption
$\gamma>4/\delta$ we have $1/\gamma<\delta/4<\delta$ so the
second term dominates the first term. The order of the variance
is thus $n^{-1+\delta}(\log n)^{1+q}$. Of course, the order of
the bias is logarithmic, hence the bias dominates the variance and the mean
squared error of $f_{nh}(x)$ is also logarithmic.
\end{rem}

\begin{rem}
Better bounds on the asymptotic variance can be obtained under
stronger mixing conditions. Consider for instance  {\em uniform
mixing}. In this case the mixing coefficient $\phi(t)$ is defined
for $t >0$ as
\begin{equation}\label{eq:betamix}
\phi(t) =  \sup_{A\in {\cal F}_{-\infty}^0, B\in {\cal F}_t^{\infty}}
|P(A|B) -P(A)|
\end{equation}
and a process is called  {\em uniform
mixing} if $\phi(t)\to 0$ for $t\to \infty$. Obviously,  uniform
mixing implies strong mixing. As a matter of fact, one has the
relation
\[
\alpha(t) \leq \tfrac{1}{2}\phi(t).
\]
See Doukhan (1994) for this inequality and many other mixing
properties. If $\sigma$ is uniform mixing with coefficient
 $\phi$ satisfying
$\int_0^\infty \phi(t)^{1/2}dt<\infty$, then  the variance
bound is given by
\begin{equation}\label{contasthm:3}
\var f_{nh}(x) = O\Big({1\over
n}\,h^{2\alpha}e^{\pi/h}\Big)+O\Big({1\over {nh\Delta}}\Big).
\end{equation}
The proof of this bound runs similarly to the strong-mixing
bound. The essential difference is that in equation~(\ref{deo})
we use Theorem 17.2.3 of Ibragimov and Linnik (1971) with
$\tau=0$  instead of Deo's (1973) lemma, as in the proof of
Theorem 2 in Masry (1983).
\end{rem}

\begin{rem}
Smoothness conditions on the density at each time of the solution
of a stochastic differential equation are guaranteed under
H\"ormander's condition, see Theorem 2.3.3 in Nualart (1995).
{Recall also relation (\ref{eq: piet})}, which can be used to
relate the smoothness of the invariant density
to the smoothness of the  drift and diffusion coefficients.
\end{rem}

\section{Proof of Theorem \ref{contasthm}}

We give the proof under the additional assumption that $b_t = 0$.
The general case is an easy consequence.
Let ${\cal F}_{\sigma}$ denote the sigma field generated by the process
$\sigma$ and let $\tilde f_{nh}$ denote the estimator based on the
approximating random variables  $\sigma_{(j-1)\Delta}Z^\Delta_j$, written as $\tilde X_j$, i.e.
\begin{equation}
\tilde f_{nh}(x)={1\over nh}\sum_{j=1}^n v_h\Big({{x- \log(
\tilde X_j^2)}\over h}\Big).
\end{equation}
The proof of (\ref{contasthm:1}) follows from the following two
lemmas, whose proofs are given in the next section. The first one
deals with the expectation  of $\tilde f_{nh}$.

\bigskip

\begin{lem}\label{lem:1}
We have
\begin{equation}\label{eq:expf}
\ex \tilde f_{nh}(x) = {1\over h}\infint w\Big({{x-u}\over h}\Big)f(u)du.
\end{equation}
\end{lem}

The second lemma estimates the expected difference between $f_{nh}$ and $\tilde f_{nh}$.
The bound is in terms of the functions
\begin{equation}\label{eq: g0}
\gamma_0(h)={1\over 2\pi }\int_{-1}^1\Big|{{\phi_w(s)
}\over\phi_k(s/h)} \Big|ds
\end{equation}
and
\begin{equation}\label{eq: g1}
\gamma_1(h,x)=e^{\frac{1}{2}\pi/h} +
{1\over h}\exp\Big(\frac{\pi}{2}\frac{1+\pi/|x|}{h}\Big)\log\frac{1+\pi/|x|}{h}.
\end{equation}

\begin{lem}\label{expbound}
For $h\to 0$ and $\eps$ small enough we have
\begin{align*}
\lefteqn{|\ex f_{nh}(x) - \ex \tilde f_{nh}(x)| =} \nonumber\\
&& O\left( {1\over h^2}\, \gamma_0(h){{\Delta}^{1/4}\over\eps} +
{1\over h}\,\gamma_0(h)\,{{\Delta^{1/2}}\over\eps^2} +
\gamma_1(h,|\log 2\eps|/h){\eps\over|\log 2\eps|}\right).
\end{align*}
\end{lem}

\bigskip

Notice that the equality~(\ref{eq:expf}) is the same as
for ordinary kernel estimators, see for instance Wand and Jones (1995). Statement
(\ref{contasthm:1}) of the theorem then follows by combining
standard arguments of kernel density estimation and
Lemma~\ref{expbound}. We will show that the bound in
Lemma~\ref{expbound} is essentially a negative power of $n$, whereas  $h^2$
is of logarithmic order. Recall that we have assumed $\delta >
4/\gamma$. It follows that ${1}/{2\gamma} < {\delta}/{4}-{1}/{2\gamma}$, so
we can pick a $\beta \in ({1}/{2\gamma}, {\delta}/{4}-{1}/{2\gamma})$ and
take $\eps=n^{-\beta}$. Up to factors that are logarithmic in
$n$ the order of $|\ex f_{nh}(x) - \ex \tilde f_{nh}(x)|$ is then
\begin{equation}
n^{\frac{1}{2\gamma}-\frac{1}{4}\delta+\beta} +
n^{\frac{1}{2\gamma}+2\beta -\frac{\delta}{2}}
+n^{\frac{1}{2\gamma}-\beta},
\end{equation}
which is negligible to $h^2=\gamma^2 \pi^2/(\log n)^2$ for the
chosen values of the parameters.
\medskip

To prove the bound (\ref{contasthm:2}) we
 use the two lemmas below, which are proved in the next section. First  consider the
variance of $\tilde f_{nh}(x)$.
\begin{lem}\label{lem:var1} We have, for $h \to 0$,
\begin{equation}
\var \tilde f_{nh}(x) = O\Big({1\over
n}\,h^{2\alpha}e^{\pi/h}\Big)+O\Big({1\over
{nh^{1+q}\Delta}}\Big).
\end{equation}
\end{lem}
The next lemma estimates $\var (f_{nh}(x)-\tilde
f_{nh}(x))$.

\begin{lem}\label{varbound}
We have, for $h\to 0$ and $\eps > 0$ small enough,
\begin{eqnarray}
\lefteqn{\var (f_{nh}(x) -   \tilde f_{nh}(x)) }\nonumber\\
&=& O\Big({1\over nh^4}\,\gamma_0(h)^2 {\Delta^{1/2}\over\eps^2}+
{1\over n}\,\gamma_1(h,|\log 2\eps|/h)^2{\eps\over|\log
2\eps|^2}\Big)\label{vbound1a}\\
&+& {1\over{nh^2\Delta}}\, O\Big({\Delta^{(1-q)/2}\over h^{2}\eps^2}\,
+\eps^{1-q}\Big).\label{vbound2a}
\end{eqnarray}
\end{lem}

The proof of (\ref{contasthm:2}) is finished as soon as we show
that the  estimate in Lemma~\ref{varbound} is of lower order than
the one in Lemma~\ref{lem:var1}. Up to terms that are
logarithmic in $n$, the bound in Lemma~\ref{lem:var1} is of order
$n^{\delta-1}$.
Choosing again $\eps=n^{-\beta}$, up to logarithmic factors, the
order of $\var( f_{nh}(x) -  \tilde f_{nh}(x))$ is
\begin{equation}
n^{-1 +\frac{1}{\gamma}-\frac{\delta}{2}+2\beta} +n^{-1
+\frac{1}{\gamma}-\beta} +n^{-1+2\beta+\frac{1+q}{2}\delta} +
n^{-1+\delta-\beta(1-q)}.
\end{equation}
Recall our assumption $\delta\gamma > 4$. If we pick $\beta$ less than $\frac{1}{4}\delta(1-q)$, then
all these terms are indeed of
lower order than $n^{\delta-1}$. \hfill$\square$

\section{Technical lemmas}
\label{sec:tech}

\subsection{Analytic properties}

We need expansions and order estimates for the functions $\phi_k$,
the kernel $v_h$, as defined in~(\ref{fourkernel}),
$\gamma_0$, as defined in (\ref{eq: g0}) and  the  function $\gamma_1$,
as defined in  (\ref{eq: g1}).
These are collected in the
lemmas of this subsection.
\begin{lem}\label{phiexpan}
 For $|t|\to\infty$ we have
\[
|\phi_k(t)|\  =\sqrt{2}\,e^{-\frac{1}{2}\pi
|t|}(1+O(\tfrac{1}{|t|})).
\]
\end{lem}

\begin{prf}
The characteristic function of $k$ is given by
\begin{equation}\label{eq: ster}
\phi_k(t)={1\over \sqrt{\pi}}\,2^{it}\,\Gamma(\tfrac{1}{2}+it).
\end{equation}
The result follows by applying the Stirling formula for the complex
gamma function, cf.
Abramowitz and Stegun (1964) Chapter 6.
\end{prf}

\begin{lem}\label{l2exp} We have the following order estimate for the
$L^2$ norm of $v_h$. For $h\to 0$
\begin{equation}
\|v_h\|_2= O(h^{\half+\alpha}e^{\pi/2h}).
\end{equation}
\end {lem}

\begin{prf}
By Parseval's identity
$$
\|v_h\|_2^2= {1\over 2\pi}\int_{-1}^1
\big|{{\phi_w(s)}\over\phi_k(s/h)}\Big|^2ds.
$$
The integral on the right-hand side is bounded  by
\begin{equation}\label{eq: abc}
\tfrac{1}{2}\int_{-1}^1|\phi_w(s)|^2 \,e^{\pi
|s/h|}ds + \int_{-1}^1|\phi_w(s)|^2\Big|{1\over|\phi_k(s/h)|^2}-
\tfrac{1}{2}e^{ \pi |s/h|}\Big|ds
\end{equation}
The first term in (\ref{eq: abc}) can be rewritten as
$$
 e^{ \pi/h}h^{1+2\alpha}\int_{0}^{1/h}\Big|{\phi_w(1-hv)
\over  (hv)^{\alpha}}
\Big|^2v^{2\alpha} \,e^{- \pi v}dv
\sim
e^{ \pi/h}h^{1+2\alpha}A^2\int_0^{\infty}v^{2\alpha} e^{- \pi v}dv,
$$
by the dominated convergence theorem. We can
rewrite the second term in   (\ref{eq: abc}) as
$$
2h^{1+2\alpha}e^{\pi/h}\int_0^{1/h}\Big|{{|\phi_w(1-hv)|}
\over (hv)^\alpha}\Big|^2 \Big|{ 2 e^{- \pi
(1/h-v)}\over|\phi_k(1/h-v)|^2}- 1\Big|v^{2\alpha} e^{- \pi
v}dv,
$$
which is of order $O(h^{1+2\alpha}e^{\pi/h})$
by the dominated convergence theorem. We have used the fact that
both the functions $\phi_w(1-u)/u^\alpha$ and $|(2\exp(-\pi
u)/|\phi_k(u)|^2)- 1|$ are bounded and that the second function
is of order $O(1/u)$ as $u$ tends to infinity. This shows that
the second  term   (\ref{eq: abc}) is negligible with respect to
the frist one.
\end{prf}
\begin{lem}\label{gamma0bound} For $h\to 0$ we have
\begin{equation}
\gamma_0(h)=O\Big(h^{1+\alpha}e^{\frac{1}{2}\pi /h}\Big).
\end{equation}
\end{lem}

\begin{prf} The proof is similar to that of Lemma~\ref{l2exp}.
\end{prf}

\bigskip

\begin{lem}\label{secondlem} The functions $v_h$ are bounded and
Lipschitz. More precisely, for all $x$  we have $|v_h(x)|\leq
\gamma_0(h)$ and for all $x$ and $u$
\begin{equation}\label{abscont}
|v_h(x+u)-v_h(x)|\leq \gamma_0(h)\, |u|.
\end{equation}
\end{lem}

\begin{prf}
The bound for $|v_h(x)|$ is obvious. To prove (\ref{abscont})
write
\[
{|v_h(x+u)-v_h(x)|}\le {1\over 2\pi}\int_{-1}^1
\Big|{\phi_w(s)\over\phi_k(s/h)}\Big|
|e^{-isu} - 1|\,ds
\le \gamma_0(h)|u|.
\]
\end{prf}

\begin{lem}\label{thirdlem} For $x\to \infty$ we have the following
estimate on the behavior of $v_h$. For some positive constant $D$ it holds that
\begin{equation}
|v_h(x)|\leq {D\over |x|}\gamma_1(h,x)\ \text{ as $ |x|\to\infty$},
\end{equation}
 and
\begin{equation}\label{gamma1bound}
\gamma_1(h,x)=O\Big({|\log h|\over h}\,e^{\frac{1}{2}\pi
(1+\pi/|x|)/h}\Big)\text{ as $h \to 0$.}
\end{equation}

\end{lem}

\begin{prf}
By a bound  in the proof of the Riemann Lebesgue lemma on page
402 of  Hewitt and Stromberg (1965) we have, with $y=\pi/x$,
\begin{eqnarray}
\lefteqn{|v_h(x)|={1\over 2\pi}\left|\infint{\phi_w(s)\over\phi_k(s/h)}\, e^{-isx}ds\right|}\nonumber\\
&\leq& {1\over 2\pi}\infint \Big|{\phi_w(s)\over\phi_k(s/h)} -
{\phi_w(s+y)\over\phi_k((s+y)/h)}\Big|\,ds\nonumber\\
&\leq& {1\over 2\pi}\infint
\Big|{{\phi_w(s)-\phi_w(s+y)}\over\phi_k(s/h)}\Big|\,ds
\label{intone}\\
&& \quad+ {1\over 2\pi}\infint \phi_w(s+y)
\Big|{1\over\phi_k(s/h)} - {1\over\phi_k((s+y)/h)}\Big|\,ds.
\label{inttwo}
\end{eqnarray}
First we need a bound on the integral (\ref{intone}). Since it
follows from Conditions W that $\phi_w$ is Lipschitz (the proof
is similar to that of~(\ref{abscont})), with Lipschitz constant
$C_1$ say, we have
\begin{eqnarray*}
\lefteqn{\infint
\Big|{{\phi_w(s)-\phi_w(s+y)}\over\phi_k(s/h)}\Big|\,ds \leq
C_1\int_{-1}^1{1\over|\phi_k(s/h)|}\, ds\, |y|}\\
&\leq& 2C_1{1\over|\phi_k(1/h)|}\, |y| \sim {C_1\over
\sqrt{2}}\,e^{\frac{1}{2}\pi/h}\, |y|,
\end{eqnarray*}
by Lemma \ref{phiexpan}. To bound the integral (\ref{inttwo}) we
need an estimate on the behaviour of $|\phi_k'|/|\phi_k|^2$.
Recall the expression (\ref{eq: ster}) for $\phi_k$.
Hence, with $\Psi=\Gamma'/\Gamma$ the digamma function,
\begin{eqnarray*}
\lefteqn{ |\phi_k'(t)|={1\over \sqrt{\pi}}\,\Big|i\log 2e^{it\log
2}\,\Gamma(\tfrac{1}{2}+it)+
i e^{it\log 2}\, \Gamma'(\tfrac{1}{2}+it)\Big|}\\
&\leq& {1\over \sqrt{\pi}}\,\Big(\log 2|\Gamma(\tfrac{1}{2}+it)|+
|\Gamma'(\tfrac{1}{2}+it)|\Big)
\end{eqnarray*}
and, as $|t|\to\infty$,
\begin{equation}\label{eq:psipsi}
\Big|{\phi_k'(t)\over\phi_k(t)^2}\Big| \leq
\sqrt{\pi}\,{1\over|\Gamma(\tfrac{1}{2}+it)|} \Big(\log 2+
|\Psi(\frac{1}{2}+it)|\Big) \leq 4\sqrt{\pi}\log(|t|)
e^{\frac{1}{2}\pi|t|},
\end{equation}
by Lemma \ref{phiexpan} and by the expansion $|\Psi(z)|\sim \log
z$ for $z\to\infty,\ |\arg z|<\pi$, cf. Abramowitz and Stegun
(1964), Chapter 6.
We now turn back to the integral~(\ref{inttwo}) and write
\begin{align*}
\infint &\phi_w(s+y) \Big|{1\over\phi_k(s/h)} -
{1\over\phi_k((s+y)/h)}\Big|\,ds\\
&= \int_{-1}^1 \phi_w(s) \Big|{1\over\phi_k((s-y)/h)} -
{1\over\phi_k(s/h)}\Big|\,ds\\
&\le {2\over h}\sup_{(-1-|y|)/h\leq s\leq (1+|y|)/h}
\Big|{\phi_k'(s)\over\phi_k(s)^2}\Big|\,
|y|\\
& \le {2\over h}\sup_{(-1-|y|)/h\leq s\leq (1+|y|)/h}
4\sqrt{\pi}\log(|s|) e^{\frac{1}{2}\pi|s|}|y| \\
&= {8\over h} \sqrt{\pi}\log((1+|y|)/h)
e^{\frac{1}{2}\pi(1+|y|)/h}|y|
\end{align*}
in view of~(\ref{eq:psipsi}). This
completes the proof.
\end{prf}

\subsection{Proof of  lemmas \ref{lem:1}-\ref{varbound}}

We start with the proof of Lemma~\ref{lem:1}. Recall that ${\cal
F}_{\sigma}$ is the $\sigma$-algebra generated by the process
$\sigma$.
\bigskip\\
{\bf Proof of Lemma~\ref{lem:1}.} Write
\begin{eqnarray*}
\lefteqn{\ex (\tilde f_{nh}(x)|{\cal F}_{\sigma})=  {1\over nh}\sum_{t=1}^n \ex
\Big(v_h\Big({{x-\log \sigma^2_{(t-1)\Delta} - \log (Z^\Delta_t)^2 }\over
h}\Big)|{\cal F}_{\sigma}\Big)}\\ &=& {1\over nh}\sum_{t=1}^n
{1\over 2\pi}\infint {{\phi_w(s)}\over\phi_k(s/h)}\ \ex
\Big(e^{-is (x-\log \sigma^2_{(t-1)\Delta} - \log (Z^\Delta_t)^2)/
h}|{\cal F}_{\sigma}\Big)ds\\ & =& {1\over nh}\sum_{t=1}^n {1\over
2\pi}\infint {{\phi_w(s)}\over\phi_k(s/h)}\, e^{-isx/h}e^{is\log
\sigma^2_{(t-1)\Delta}/h}\,\phi_k(s/h)ds\\ &=& {1\over
nh}\sum_{t=1}^n{1\over 2\pi}\infint \phi_w(s)\ e^{-is(x-\log
\sigma^2_{(t-1)\Delta})/h}ds
\\
&=& {1\over nh}\sum_{t=1}^n w\Big({{x-\log
\sigma^2_{(t-1)\Delta}}\over h}\Big).
\end{eqnarray*}
By taking expectation the result follows. \hfill$\square$
\medskip

For the proof of Lemma~\ref{expbound} we need
a few properties of the process $\sigma$, valid under Condition
$\sigma$. Since $(x-y)^2 \le |x^2-y^2|$ for $x,y \ge 0$, it holds  that
$\ex(\sigma_t-\sigma_0)^2 =O(t^{1/2})$
for $t \to 0$. Consequently, there exists a constant $C> 0$
such that
\begin{equation}\label{sigma2}
 \ex (X_1^\Delta -\sigma_{0}Z^\Delta_1)^2 \leq C \Delta^{1/2} \text{ for $\Delta \to 0$},
\end{equation}
since $\ex (X_1^\Delta -\sigma_{0}Z^\Delta_1)^2=\frac{1}{\Delta}\ex
\int_0^\Delta (\sigma_t-\sigma_0)^2\,dt$. Moreover, Condition $\sigma$ implies that
\begin{equation}\label{sigma3}
\ex \left|\tfrac{1}{\Delta}\int_0^\Delta \sigma^2_t\,dt -
\sigma_0^2 \right| =O(\Delta^{1/2})\text{ for $\Delta \to 0$}.
\end{equation}
\medskip \\
{\bf Proof of Lemma~\ref{expbound}.} Writing
$$
W_j=v_h\Big({{x-\log((X_j^\Delta)^2)}\over h}\Big)-
v_h\Big({{x-\log(\tilde X_j^2)}\over h}\Big),
$$
so that $f_{nh}(x) -   \tilde f_{nh}(x)={1\over {nh}}\sum_{j=1}^n W_j$,
we have
\begin{eqnarray}
\lefteqn{|\ex f_{nh}(x) - \ex \tilde f_{nh}(x)|
\leq
 {1\over h}\,\ex |W_j|}\nonumber\\
&=& {1\over h}\,\ex |W_j|
 I_{[|X_1^\Delta|\geq\eps\ \mbox{and}\
|\tilde X_1|\geq\eps]}\label{partone}\\
&\quad& +{1\over h}\,\ex |W_j|
 I_{[|X_1^\Delta|\leq\eps\ \mbox{or}\
|\tilde X_1|\leq\eps]} I_{[|X_1^\Delta -\tilde X_1|\geq\eps]}
\label{parttwo}\\
&\quad& +{1\over h}\,\ex |W_j|
 I_{[|X_1^\Delta|\leq\eps\ \mbox{or}\
|\tilde X_1|\leq\eps]} I_{[|X_1^\Delta -\tilde X_1|<\eps]}.
\label{partthree}
\end{eqnarray}
By Lemma \ref{secondlem} and (\ref{sigma2}) the term
(\ref{partone}) can be bounded by
\begin{eqnarray*}
\lefteqn{{2\over h^2}\,\gamma_0(h)\ex
|\log((X_1^\Delta)-\log(\tilde X_1)| I_{[|X_1^\Delta|\geq\eps\
\mbox{and}\ |\tilde
X_1|\geq\eps]}}\\
&\leq& {2\over h^2}\,{1\over\eps}
\gamma_0(h)\ex |X_1^\Delta-\tilde X_1|
\leq {2\over h^2}\,
\gamma_0(h)\sqrt{C}{{\Delta}^{1/4}\over\eps}.
\end{eqnarray*}
In the same way the term (\ref{parttwo}) can be bounded by
$$
{2\over h}\,\gamma_0(h)P(|X_1^\Delta -\tilde X_1|\geq\eps)
\leq {2\over h}\,\gamma_0(h)C\,{\Delta^{1/2}\over\eps^2}.
$$
Since the absolute value of both arguments of $v_h$ below are
 eventually larger than $|\log 2\eps |/h$,
by Lemma \ref{thirdlem} the term (\ref{partthree}) can be bounded
by
$$
{1\over h}\,\gamma_1(h,|\log 2\eps|/h){1\over(|\log
2\eps|/h)}\,P(|\tilde X_1|\leq
2\eps)
\leq C_2 \,\gamma_1(h,|\log 2\eps|/h){\eps\over|\log 2\eps|},
$$
for some constant $C_2$. Here we used the fact that the density of $\tilde{X}_1$ is
bounded which follows from the assumption that
 $\sigma^2_0$ has a bounded density in a
neighbourhood of zero.
\hfill$\square$

\bigskip

\noindent {\bf Proof of Lemma~\ref{lem:var1}.} Consider the
decomposition
\begin{equation}\label{vardecomp}
\var (\tilde f_{nh}(x))= \var(\ex (\tilde f_{nh}(x)|{\cal
F}_{\sigma})) + \ex(\var (\tilde f_{nh}(x)|{\cal F}_{\sigma})).
\end{equation}
By the proof of Lemma \ref{lem:1} the conditional expectation $\ex (\tilde f_{nh}(x)|{\cal
F}_{\sigma})$ is equal to a
kernel estimator of the density of $\log \sigma^2_{t}$. By
Theorem 3 of Masry (1983), we can bound its variance by
$$
 {{20(1+o(1))}\over {nh^{1+q}\Delta}}\ f(x)^{1-q}
\Big(\infint|w(u)|^{2/(1-q)}du \Big)^{1-q} \int_0^\infty
\alpha(\tau)^qd\tau=O\Big({1\over {nh^{1+q}\Delta}}\Big).
$$
Given the process $\sigma$
 the
random variables $\log \tilde{X}^2_t$ are independent, so we can
bound  the second term in (\ref{vardecomp})  by
$$
{1\over n^2h^2}\sum_{t=1}^n \var \Big(v_h\Big({{x-\log
\tilde{X}^2_t}\over h}\Big)\Big)
\leq {1\over nh^2}\, \ex \Big(v_h\Big({{x-\log
\tilde{X}^2_1}\over h}\Big)\Big)^2\leq {1\over nh^2}\gamma_0(h)^2,
$$
by Lemma \ref{secondlem}.
The result follows by an application of Lemma~\ref{gamma0bound}.
\hfill$\square$
\medskip\\
{\bf Proof of Lemma~\ref{varbound}.}
Note that for different $i,j$, conditional on the process
$\sigma$,
 the pairs
$X_i^\Delta,\tilde X_i$ and $X_j^\Delta,\tilde X_j$ are
independent. Hence the conditional covariances of functions of
these pairs vanish.

With $W_j$ as in the proof of Lemma \ref{expbound} we have
\begin{eqnarray}
\lefteqn{\var( f_{nh}(x) -   \tilde f_{nh}(x))}\nonumber\\
&=& {1\over nh^2}\ \var W_1 +
 {1\over n^2h^2}\sum_{i\not= j} \cov(\ex( W_i|{\cal F}_\sigma),\ex(
W_j|{\cal F}_\sigma))\label{sumbound}.
\end{eqnarray}
Let us first derive a bound on $\Var W_1$. We have $\Var W_1 \le \ex
W_1^2$,
which can be split up in three terms
\begin{eqnarray}
 &\quad& {1\over h}\,\ex W_j^2
 I_{[|X_1^\Delta|\geq\eps\ \mbox{and}\
|\tilde X_1|\geq\eps]}\label{varpartone}\\
&\quad& +{1\over h}\,\ex W_j^2
 I_{[|X_1^\Delta|\leq\eps\ \mbox{or}\
|\tilde X_1|\leq\eps]} I_{[|X_1^\Delta -\tilde X_1|\geq\eps]}
\label{varparttwo}\\
&\quad& +{1\over h}\,\ex W_j^2
 I_{[|X_1^\Delta|\leq\eps\ \mbox{or}\
|\tilde X_1|\leq\eps]} I_{[|X_1^\Delta -\tilde X_1|<\eps]}.
\label{varpartthree}
\end{eqnarray}
By (\ref{sigma2}) and Lemma \ref{secondlem} the term
(\ref{varpartone}) can be bounded by
${2\over h^2}\, \gamma_0(h)^2\,C\,{ \Delta^{1/2}
\over\eps^2}$.
Again by (\ref{sigma2}) and Lemma \ref{secondlem} the term
(\ref{varparttwo}) can be bounded by
\[
4\gamma_0(h)^2P(|X_1^\Delta -\tilde X_1|\geq\eps)
\le 4\gamma_0(h)^2C\,{\Delta^{1/2}\over\eps^2}.
\]
Since the absolute value of both arguments of $v_h$ below are
 eventually larger than $|\log 2\eps |/h$,
by Lemma \ref{thirdlem} the term (\ref{varpartthree}) can be
bounded by
$$
{{\gamma_1(h,|\log 2\eps|/h)^2 }\over{(|\log
2\eps|/h)^2}}\,P(|\tilde X_1|\leq
2\eps)
\leq C_2 \,h^2\gamma_1(h,|\log 2\eps|/h)^2{\eps\over|\log
2\eps|^2},
$$
for some constant $C_2$, where we use  again, as in the  proof of Lemma~\ref{expbound}, that the density of
$\tilde{X}_1$ is bounded.
\\
We get
\begin{equation}\label{expz}
\ex W_1^2 = O\Big(  {1\over h^2}\, \gamma_0(h)^2\,C\,{
\Delta^{1/2} \over\eps^2} +h^2\gamma_1(h,|\log
2\eps|/h)^2{\eps\over|\log 2\eps|^2}\Big),
\end{equation}
which gives the first order bound (\ref{vbound1a}).

Next we concentrate on the sum of covariances in (\ref{sumbound}).
Define
\begin{equation}
\bar{\sigma_i} = {1\over\Delta}\int_{(i-1)\Delta}^{i\Delta}
\sigma_t^2dt.
\end{equation}
Note that given ${\cal F}_\sigma$, $X_i^\Delta$ is ${\cal
N}(0,\bar{\sigma_i})$ distributed and $\tilde{X_i}$ is
$N(0,\sigma^2_{(i-1)\Delta})$.  As in the proof of Lemma~\ref{lem:1} it
follows that
$$
\ex( W_i|{\cal F}_\sigma)=w\Big({{x-\log \bar\sigma_i}\over
h}\Big) -w\Big({{x-\log \sigma_{(i-1)\Delta}^2}\over h}\Big).
$$

We follow the line of arguments in the proof of Theorem 3 in
Masry (1983). The stationarity of $W_j$ implies that also the
conditional expectations $\tilde{W}_j:=\ex(W_j|{\cal F}_\sigma)$
are stationary. Hence we have
\[
\sum_{i\not= j}\cov(\tilde{W}_i,\tilde{W}_j)=
2 \sum_{k=1}^{n-1}(n-k)\cov(\tilde{W}_0,\tilde{W}_k).
\]
Now note that the process $\tilde{W}_j$ is strongly mixing with a
mixing coefficient $\tilde\alpha (k)\leq \alpha( (k-1)\Delta),
k=1,2,\dots$, where $\alpha$ is the coefficient of the process
$\sigma$. By a lemma of Deo (1973) for strongly mixing processes
it follows that for all $\tau>0$
\begin{equation}\label{deo}
|\cov(\tilde{W}_0,\tilde{W}_k)| \le 10\alpha((k-1)\Delta)^{\tau/(2+\tau)}
\Big(\ex|\tilde{W}_1|^{2+\tau}\Big)^{{2/(2+\tau)}}.
\end{equation}
By the monotonicity of the mixing coefficient $\alpha$ we get
\begin{eqnarray*}
\lefteqn{\Big|{1\over n^2h^2}\sum_{i\not= j} \cov(\tilde{W}_i,
\tilde{W}_j)\Big|}\\
&\leq& {10\over
nh^2}\,\Big(\ex|\tilde{W}_1|^{2+\tau}\Big)^{{2/(2+\tau)}}
 \sum_{k=1}^{n-1}(1-{k\over
n})\alpha((k-1)\Delta)^{\tau/(2+\tau)}\\
&\leq& {10\over nh^2}\Big(\alpha(0)^{\tau/(2+\tau)}+{1\over
\Delta}\int_0^\infty\alpha(t)^{\tau/(2+\tau)}dt\Big)
\Big(\ex|\tilde{W}_1|^{2+\tau}\Big)^{{2/(2+\tau)}}.
\end{eqnarray*}
Next we derive a bound on  $\ex|\tilde{W}_1|^{2+\tau}$. Fix $\kappa\in (0,1]$. We have
\begin{eqnarray}
\lefteqn{\ex|\tilde{W}_1|^{2+\tau}
= \ex \Big|w\Big({{x-\log(\bar\sigma_1)}\over h}\Big)-
w\Big({{x-\log(\sigma_{0}^2)}\over h}\Big)\Big|^{2+\tau}\nonumber}\\
&&\quad\quad\quad\times I_{[\bar\sigma_1\geq\eps\ \mbox{and}\
\sigma_{0}^2\geq\eps]}
\label{varwpartone}\\
&+& \ex \Big|w\Big({{x-\log(\bar\sigma_1)}\over h}\Big)-
w\Big({{x-\log(\sigma_{0}^2)}\over h}\Big)\Big|^{2+\tau}\nonumber\\
&&\quad\quad\quad\times I_{[\bar\sigma_1\leq\eps\ \mbox{or}\
\sigma_{0}^2\leq\eps]} I_{[|\bar\sigma_1^\kappa -\sigma_{0}^{2\kappa}|\geq\eps]}
\label{varwparttwo}\\
&+& \ex \Big|w\Big({{x-\log(\bar\sigma_1)}\over h}\Big)-
w\Big({{x-\log(\sigma_{0}^2)}\over h}\Big)\Big|^{2+\tau}\nonumber\\
&&\quad\quad\quad\times I_{[\bar\sigma_1\leq\eps\ \mbox{or}\
\sigma_{0}^2\leq\eps]} I_{[|\bar\sigma_1^\kappa -\sigma_{0}^{2\kappa}|<\eps]}.
\label{varwpartthree}
\end{eqnarray}
Note that by Condition W and Fourier inversion  $w$ is Lipschitz with constant
$1/\pi$ and bounded by $1/\pi$. Hence the term (\ref{varwpartone}) can be bounded by
$\ex|\bar\sigma_1^\kappa-\sigma_{0}^{2\kappa}|^{2+\tau}/(\kappa\eps h)^{2+\tau}$.
The term
(\ref{varwparttwo}) can be bounded by
\begin{equation*}
 P(|\bar\sigma_1^\kappa -\sigma_{0}^{2\kappa}|\geq\eps)
\leq {1\over\eps^{2+\tau}}\,\ex
|\bar\sigma_1^\kappa-\sigma_{0}^{2\kappa}|^{2+\tau}.
\end{equation*}
Likewise, the term (\ref{varwpartthree}) can be bounded by
$$
\ex \Big|w\Big({{x-\log(\bar\sigma_1)}\over h}\Big)-
w\Big({{x-\log(\sigma_{0}^2)}\over h}\Big)\Big|^{2+\tau}
I_{[\bar\sigma_1\leq \eps(1+\eps^{1-\kappa})^{1/\kappa}\ \mbox{and}\
\sigma_{0}^2\leq
\eps(1+\eps^{1-\kappa})^{1/\kappa}]},
$$
which is bounded by $P(\sigma_{0}^2\leq
2\eps)=O(\eps)$ since $\sigma^2_0$ was assumed to have a bounded density in a
neighbourhood of zero.
\\
With $\tau=2q/(1-q)$ and $\kappa=\frac{1}{2+\tau}=\frac{1-q}{2}$
 we have  with an application of the basic inequality
$|u^\kappa-v^\kappa|\leq |u-v|^\kappa$ for $u, v \geq 0$ and $\kappa\in(0,1]$ in
the second equality below and from condition $\sigma$ and its
consequence~(\ref{sigma3}) in the fourth equality
\begin{eqnarray*}
\lefteqn{\Big|{1\over n^2h^2}\sum_{i\not= j} \cov(\tilde{W}_1,
\tilde{W}_j)\Big|}\\
&=& \frac{1}{nh^2\Delta}O\Big({1\over
h^{2+\tau}}\,{1\over\eps^{2+\tau}}\, \ex
|\bar\sigma_1^\kappa-\sigma_{0}^{2\kappa}|^{2+\tau}
+\eps\Big)^{2/(2+\tau)}\\
&=&  \frac{1}{nh^2\Delta}O\Big({1\over
h^{2+\tau}}\,{1\over\eps^{2+\tau}}\, \ex
|\bar\sigma_1-\sigma_{0}^{2}|^{\kappa(2+\tau)}
+\eps\Big)^{2/(2+\tau)}  \\
&=& \frac{1}{nh^2\Delta} O\Big({(\ex
|\bar\sigma_1-\sigma_{0}^{2}|)^{2/(2+\tau)}\over {h^{2}\eps^{2}}}\,
+\eps^{2/(2+\tau)}\Big), \\
&=& \frac{1}{nh^2\Delta} O\Big({\Delta^{1/(2+\tau)}\over {h^{2}\eps^{2}}}\,
+\eps^{2/(2+\tau)}\Big) \\
&=& \frac{1}{nh^2\Delta} O\Big({\Delta^{(1-q)/2}\over {h^{2}\eps^{2}}}\,
+\eps^{1-q
}\Big),\end{eqnarray*}
which gives the second order bound (\ref{vbound2a}).
\hfill $\square$

\section*{References}

\small

\begin{verse}

Abramowitz, M.\ and Stegun, I\ (1964) , {\em Handbook of
Mathematical Functions, ninth edition}, Dover, New York.






Deo, C.M.\ (1973), A note on empirical processes for strong
mixing processes, {\em Ann. Probab.} {\bf 1}, 870--875.

Doukhan, P.\ (1994), {\em Mixing, Properties and Examples},
Springer-Verlag.

Fan, J.\ (1991), On the optimal rates of convergence for
nonparametric deconvolution problems, {\em Ann. Statist.} {\bf
19}, 1257--1272.


Genon-Catalot, V., Jeantheau, T.\ and Lar\'edo, C.\ (1998), Limit
theorems for discretely observed stochastic volatility models,
{\em Bernoulli} {\bf 4}, 283--303.

Genon-Catalot, V., Jeantheau, T.\ and Lar\'edo, C.\ (1999),
Parameter estimation
 for discretely observed stochastic volatility models,
{\em Bernoulli} {\bf 5}, 855-872.

Genon-Catalot, V., Jeantheau, T.\ and Lar\'edo, C.\ (2000),
Stochastic volatility models as hidden Markov models and
statistical applications, {\em Bernoulli} {\bf 6}, 1051--1079.

Gihman, I.I. and  Skorohod A.V.\ (1972), {\em Stochastic
Differential Equations}, Springer.

Heston, S.L.\  (1993), A closed-form solution for options with
stochastic volatility with applications to Bond and Currency
options, {\em The Review of Finacial Studies} {\bf 6} (2),
327--343.

Hewitt, E.\ and Stromberg K.\ (1965), {\em Real and Abstract
Analysis}, Springer Verlag, New York.

Ibragimov, I.A., and Linnik (1971), {\em Independent and
stationary sequences of random variables}, Wolters-Noordhoff.

Karatzas, I.\  and S.E.\ Shreve (1991), {\em Brownian Motion and
Stochastic Calculus}, Springer Verlag, New York.

Masry, E.\ (1983), Probability density estimation from sampled
data, {\em IEEE Trans. Inform. Theory}  {\bf 29}, 696--709.



Masry, E.\ (1993), Strong consistency and rates for
deconvolution of multivariate densities of stationary processes,
{\em Stoc. Proc. and Appl.} {\bf 475}, 53--74.

Nualart, D.\ (1995), {\em The Malliavin calculus and related
topics}, Springer Verlag, New York.


Skorokhod, A.V.\ (1989), {\em Asymptotic Methods in the Theory of
Stochastic Differential Equations}, AMS.


Wand, M.P.\ (1998), Finite sample performance of deconvolving
kernel density estimators, {\em Statist.\ Probab.\ Lett.} {\bf 37},
131--139.

Wand, M.P.\ and Jones, M.C.\ (1995), Kernel Smoothing, Chapman and Hall,
London.

Wiggins, J.\ B.\  (1987), Option valuation under stochastic
volatility, {\em Journal of Financial Economics} {\bf 19},
351--372.

\end{verse}

\end{document}